\documentclass[11pt]{article}
\usepackage[english]{babel}
\usepackage{amsmath}
\usepackage{amsthm}
\usepackage{amssymb}
\usepackage[dvips]{color}
\newtheorem{theorem}{Theorem}[section]
\newtheorem{definition}[theorem]{Definition}
\newtheorem{lemma}[theorem]{Lemma}
\newtheorem{prop}[theorem]{Proposition}
\newtheorem{corollary}[theorem]{Corollary}
\newtheorem{remark}[theorem]{Remark}
\newtheorem{example}[theorem]{Example}
\newcommand{\pn}{\par\noindent}

\newcommand{\Er}{\mathcal{R}}
\newcommand{\ov}{\overline}

\newcommand{\med}{\medskip}
\newcommand{\rr}{\mathbb{R}}

\newcommand{\hh}{\mathbb{H}}

\newcommand{\pp}{\partial}

\title{\bf    Non commutative functional calculus: bounded operators  }

\author{Fabrizio Colombo\\Dipartimento di Matematica\\Politecnico di
Milano\\Via Bonardi, 9\\20133 Milano,
Italy\\fabrizio.colombo@polimi.it \and Graziano
Gentili\\Dipartimento di Matematica\\ Universit\'a di Firenze
\\Viale Morgagni, 67 A\\
Firenze, Italy\\ gentili@math.unifi.it
\\
\and Irene Sabadini\\Dipartimento di
Matematica\\Politecnico di Milano\\Via Bonardi, 9\\20133 Milano,
Italy\\irene.sabadini@polimi.it
\and Daniele C. Struppa\\Department of Mathematics \\and Computer Sciences
\\Chapman University\\Orange, CA 92866 USA,
\\struppa@chapman.edu}

\begin{document}

\maketitle
\begin{abstract}
In this paper we  develop a functional calculus for bounded
operators defined on quaternionic Banach spaces. This calculus is
based on the notion of slice-regularity, see \cite{gs}, and
the key tools are a new resolvent operator and a new eigenvalue
problem.
\end{abstract}



\section{Introduction}

Let $V$ be a Banach space over a (possibly skew) field
$\mathbb{F}$ and let $E$ be the Banach space of linear operators
acting on it. If $T\in E$, and if $\mathbb{F}=\mathbb{C}$, then
the standard eigenvalue problem seeks the values
$\lambda\in\mathbb{C}$ for which $\lambda \mathcal{I}-T$ is not
invertible. This situation is well known and leads to what is
known as Functional Calculus (F.C.). Indeed, in the complex case,
one defines the spectrum of an operator as the set for which
 $\lambda \mathcal{I}-T$ is not invertible. The key observation for the
development of the F.C. is the fact that the inverse of
$\lambda \mathcal{I}-T$ formally coincides with a Cauchy kernel
which allows integral representation for holomorphic functions. As
a consequence, one is able to consider for any function $f$
holomorphic on the spectrum of $T$ its value $f(T)$ which is
formally defined as
$$
f(T)=\frac{1}{2\pi i}\int_\Gamma (\lambda \mathcal{I}- T)^{-1} f(\lambda)
d\lambda
$$
where $\Gamma$ is a closed curve surrounding the spectrum of $T$.
 The important fact is that such a
definition coincides with the obvious meaning when $f$ is a
polynomial and that this definition behaves well with respect to
linear combinations, product and composition of functions, i.e.,
$(af+bg)(T)=af(T)+bg(T)$, $a,b\in\mathbb{C}$, $(fg)(T)=f(T)g(T)$,
$(f\cdot g)(T)=f(g(T))$. For the basic theory in the complex case
a classical reference is \cite{ds}. Important applications of this
theory to specific operators are now an important field of
investigation (see for example the references in
\cite{jefferies}).

Several new difficulties arise when one tries to deal with
quaternionic linear operators. First of all, when working in a non
commutative setting, it is necessary to specify that the operators
are linear, for example, on the right, i.e.
$T(v\alpha+w\beta)=T(v)\alpha +T(w)\beta$, for
$\alpha,\beta\in\mathbb{H}$. In addition, it turns out that there
are two different eigenvalue problems. The so called right
eigenvalue problem, i.e. the search for non zero vectors
satisfying $T(v)=v\lambda$, is widely studied by physicists. The
crucial fact is that whenever there is a non real eigenvalue
$\lambda$ then all quaternions belonging to the sphere
$r^{-1}\lambda r$, $r\in\mathbb{H}\setminus \{0\}$, are
eigenvalues. This fact allows to choose the phase and to work, for
example, with complex eigenvalues.

Note however that the operator of multiplication on the right is
not a right linear operator, and so the operator
$\mathcal{I}\lambda -T$ is not linear. On the other hand, the
operator $\lambda\mathcal{I}-T$ is right linear but, in general,
one cannot choose the phase as the eigenvalues do not necessarily
lie on a sphere. Even more important, in the complex setting the
inverse $(\lambda\mathcal{I}-T)^{-1}$ is related to a Cauchy
kernel useful in the notion of Cauchy integral. In the
quaternionic setting it is not even clear which type of
generalized holomorphy must be used. The regularity in the sense
of Fueter (see e.g. \cite{nostro libro}) does not seem to give a
good notion of exponential function and does not allow to
introduce polynomials (or even powers) of operators. While many
authors have looked at the quaternionic eigenvalue problem because
of its applications to physics, \cite{adler}, it is safe to say
that until now there has been no general treatment which could
allow the construction of a F.C. with quaternionic spectrum, as we
describe in this paper.

As we will show, however, a recent new theory, the theory of
slice-regular functions, see \cite{gs}, \cite{advances} where they are called
Cullen-regular functions, allows to show that power series in the
quaternionic variable are regular, so this new notion seems to be
the correct setting in which a F.C. can be introduced and
developed. The Cauchy kernel series $\sum_{n\geq 0} q^ns^{-n-1}$,
$q,s\in\mathbb{H}$, which is used to write a Cauchy formula for
slice-regular functions does not coincide, in general, with
$(s-q)^{-1}$, thus the linear operator $S^{-1}(s,T)=\sum_{n\geq 0}
T^ns^{-n-1}$ is not, in general, the inverse of
$(s\mathcal{I}-T)$. The key idea is to identify the operator whose
inverse is $\sum_{n\geq 0} T^ns^{-n-1}$. This new operator will
give us a new notion of spectrum (the so called S-spectrum) which
will allow to introduce a F.C. As we will see, this spectrum, as
the right-spectrum, will allow the choice of the phase. In this
paper we confine our attention to the case in which the operators
under consideration are bounded. However, as we will show in a
subsequent paper, much of this theory can be extended to the
unbounded case (which of course has significant physical
applications).

We close this introduction by pointing out that the spectral
theory can be extended to the case of $n$-tuples of operators. For
the complex case, the reader is referred to \cite{Taylor} (for the
commuting case) and to \cite{jefferies} for the case in which the
operators do not commute. We plan to come back to this issue in a
subsequent paper and show  how our theory can be applied to this
situation as well.

{\it Acknowledgements} The first and third authors are grateful to
Chapman University for the hospitality during the period in which
this paper was written. They are also indebted to G.N.S.A.G.A. of INdAM and the
Politecnico di Milano for partially supporting their visit.

\section{Regular functions and linear bounded operators}

\subsection{\it Regular functions.}

In this section we collect some basic results that we need in the
sequel. For more details we refer the reader to [11].

\pn
Let $\hh$ be the real associative algebra of quaternions
with respect to the basis $\{1, i,j,k \}$
satisfying the relations
$$
i^2=j^2=k^2=-1,
 ij =-ji =k,\
$$
$$
jk =-kj =i ,
 \  ki =-ik =j .
$$
We will denote a quaternion $q$ as $q=x_0+ix_1+jx_2+kx_3$,
$x_i\in \mathbb{R}$, its conjugate as
$\bar q=x_0-ix_1-jx_2-kx_3$, and we will write $|q|^2=q\ov q$.
\noindent
In the recent papers \cite{gs}, \cite{advances}, the authors introduced a
new notion of regularity for functions on $\hh$. We will recall here both the definition
and its most salient consequences.
\noindent
Let $\mathbb{S}$ be the sphere of purely imaginary unit quaternions, i.e.
$$
\mathbb{S}=\{ q=ix_1+jx_2+kx_3\ |\  x_1^2+x_2^2+x_3^2=1\}.
$$
\begin{definition} Let $U\subseteq\hh$ be an open set and let
$f:\ U\to\hh$ be a real differentiable function. Let $I\in\mathbb{S}$ and let $f_I$ be
the restriction of $f$ to the complex line $L_I := \mathbb{R}+I\mathbb{R}$ passing through $1$ and $I$.
We say that $f$ is a left regular function if for every $I\in\mathbb{S}$
$$
\frac{1}{2}\left(\frac{\partial }{\partial x}+I\frac{\partial }{\partial y}\right)f_I(x+Iy)=0,
$$
and we say it is right regular if for every $I\in\mathbb{S}$
$$
\frac{1}{2}\left(\frac{\partial }{\partial x}+\frac{\partial }{\partial y}I\right)f_I(x+Iy)=0.
$$
\end{definition}

\begin{remark}{\rm Left regular functions on $U\subseteq\hh$ form a right vector space $\Er (U)$ and
right regular functions on $U\subseteq\hh$ form a left vector space.
 It is not true, in general, that the product of two
regular functions is regular.}
\end{remark}

\begin{example}{\rm The key feature of this notion of regularity is the fact that
polynomials $\sum_{n=0}^Nq^na_n$ in the quaternion variable $q$, and with quaternionic
coefficients $a_n$ are left regular (while polynomials $\sum_{n=0}^Na_nq^n$ are right regular).
Moreover any power series $\sum_{n=0}^{+\infty} q^na_n$ is left regular in
its domain of convergence. As an example, the function $R(q)=(q-q_0)^{-1}$, $q_0\in\mathbb{R}$
which corresponds to the Cauchy kernel (see below), is regular for $|q|<|q_0|$.}
\end{example}
\noindent
In fact, every regular function can be represented as a power series, \cite{advances}:
\begin{theorem} If $B=B(0,R)$ is the open ball centered in the origin with radius $R>0$ and
$f:\ B \to\hh$ is a left regular function, then $f$ has a series expansion of the form
$$
f(q)=\sum_{n=0}^{+\infty} q^n\frac{1}{n!}\frac{\pp^nf}{\pp x^n}(0)
$$
converging on $B$. Analogously, if $f$ is right regular it can be expanded as
$$
f(q)=\sum_{n=0}^{+\infty} \frac{1}{n!}\frac{\pp^nf}{\pp x^n}(0) q^n.
$$
\end{theorem}
\noindent
Note that, even though the definition of regular
function involves the direction of the unit quaternion $I$, the coefficients of the series
expansion do not depend at all from a choice of $I$.
\begin{remark}{\rm From now on we will not specify whether we are considering left or right regular
functions, since the context will clarify it.}\end{remark}
\noindent
A main result in the theory of regular functions is the analogue of the Cauchy integral formula.
In order to state the result we need some notation.
Given a quaternion $q=x_0 + ix_1 + j x_2 + kx_3$ let us denote its real part
$x_0$ by $Re[q]$ and its imaginary part $ix_1 + j x_2 + kx_3$ by $Im[q]$.
We set
$$
I_q=\left\{\begin{array}{l}
\displaystyle\frac{Im[q]}{|Im[q]|}\quad{\rm if}\ Im[q]\not=0\\
{\rm any\ element\ of\ } \mathbb{S}{\rm\ otherwise}\\
\end{array}
\right.
$$
We have the following \cite{advances}:
\begin{theorem}
Let $f:\ B(0,R)\to\hh$ be a regular function and let $q\in B(0,R)$. Then
$$
f(q)=\frac{1}{2\pi } \int_{\pp\Delta_q(0,r)}
\displaystyle (\zeta -q)^{-1}\, d\zeta_{I_q} \ f(\zeta)
$$
where $d\zeta_{I_q}=-I_qd\zeta$ and $r>0$ is such that
$$
\overline{\Delta_q(0,r)}:=\{x+I_q y\ |\ x^2+y^2\leq r^2\}
$$
contains $q$ and is contained in $B(0,R)$.
\end{theorem}
\noindent
The proof of the theorem relies on the following result, which is of independent interest :
\begin{theorem}\label{int_nullo}
Let $f:\ B(0,R)\to\hh$ be a real differentiable function and let $q\in B(0,R)$.
If $f$ is regular and $r$ is such that $\overline{\Delta_q(0,r)}$ is contained in $B(0,R)$, then
$$
\int_{\pp\Delta_q(0,r)} d\zeta f(\zeta)=0 .
$$
\end{theorem}
\pn
\subsection{\it Regular functions with values in a Banach space.}

We now turn our attention to functions $f:\ \hh\to E$ where $E$ denotes
a left quaternionic Banach space. Let us revise the definition of integral in this setting.

\begin{definition}
Let $E$ be a left quaternionic Banach space.
A function $f:\ \hh\to E$ is said to be
regular if there exist an open ball $B=B(0,R)\subseteq \hh$ and a sequence $\{a_n\}$ of elements of $E$
such that, for every point $q\in B(0,R)$ the function $f(q)$
can be represented by the following  series
\begin{equation}\label{(2.14)}
f(q)=\sum_{n=0}^{+\infty} q^n a_n,
\end{equation}
 converging in the norm of $E$ for any $q$ such that $|q|< R$.
\end{definition}
\pn
Consider $E'={\rm Hom}(E,\hh)$,  the dual space of $E$, and let  $x'\in E'$ be a
linear continuous functional. If $\Gamma$ is any compact set of $B$, then
by the left linearity and the continuity of $x'$
we have
\begin{equation}\label{2.15}
< \int_{\Gamma} d\zeta\, f(\zeta) ,x'> = \int_{\Gamma}   d\zeta\, <f(\zeta) ,x'> .
\end{equation}
To prove the  Cauchy integral formulas in the case of vector functions
we need the following theorem which is a consequence of the Hahn-Banach
theorem.
\med
\pn

\begin{theorem}\label{Hahn-Banach} Let $E$ be a left quaternionic Banach space
and let $x\in E$. If for every linear and
continuous functional $x'$ we have $<x,x'>=0$ then $x=0$.
\end{theorem}
\pn
\begin{proof}
We observe that the Hahn-Banach theorem still holds true
for left (or  right) vector spaces on $\hh$.
Let $x\in E$, $x$ non zero and let $W$ be the subspace generated by $x$.
Let $\xi':\ W\to \hh$ be the functional such that $\xi'(xq)=\| x\|q$. Obviously, $\xi'$ is left linear continuous
functional of norm $1$. By the Hahn-Banach theorem there exists an extension $x'$ of $\xi'$ to all of
$E$ with $\|x'\|=1$, moreover $x'(x)=\| x\|$ and the statement follows.
\end{proof}
\pn
\begin{theorem} Let $E$ be a left quaternionic
Banach space and let $f: \hh \to E$
be a regular function on an open ball $B(0,R)$ containing $\Delta_q(0,r)$. Then
\begin{equation}\label{2.19}
\int_{\pp\Delta_q(0,r)}  d\zeta\, f(\zeta)=0 .
\end{equation}
\end{theorem}
\begin{proof} For every $x'\in E'$ and for any $q\in B(0,R)$ we have:
$$
f(q)=\sum_{n=0}^{+\infty} q^n a_n ,\ \ \ \ \ \ \ \ \ (a_n\in E)
$$
$$
<f(q),x'>=<\sum_{n=0}^{+\infty}q^n a_n,x'>=\sum_{n=0}^{+\infty}q^n <a_n,x'> .
$$
Hence, by the Abel's lemma for quaternionic power series, the function $<f(q),x'>$ is regular in
$B(0,R)$. Thanks to  Theorem \ref{int_nullo}, we have
$$
\int_{\pp\Delta_q(0,r)}d\zeta\, < f(\zeta)  , x'>\ \ \  = 0.
$$
By the equalities
\begin{equation}\label{2.20}
0=\int_{\pp\Delta_q(0,r)} d\zeta\, < f(\zeta)  , x'>=< \int_{\pp\Delta_q(0,r)} d\zeta\,  f(\zeta) , x' >,
\ \ \ \ \forall x'\in E'
\end{equation}
and from Theorem \ref{Hahn-Banach} we obtain (\ref{2.19}).\end{proof}
\pn
\begin{theorem}
Let $E$ be a left Banach space and let $f: \hh \to E$
be a regular function on $B(0,R)$.  Then
$$
f(q)=\frac{1}{2\pi } \int_{\pp\Delta_q(0,r)}
(\zeta -q)^{-1}\, d\zeta_{I_q} \ f(\zeta )
$$
where $d\zeta_{I_q}=-I_qd\zeta$ and $r>0$ is such that
\begin{equation}\label{2.21}
\overline{\Delta_q(0,r)}=\{x+I_q y\ |\ x^2+y^2\leq r^2\}
\end{equation}
contains $q$ and is contained in $B(0,R)$.
\end{theorem}
\begin{proof}
The statement follows as in the complex case.
\end{proof}

\subsection{\it Linear bounded quaternionic operators.}
We conclude this section with a quick discussion of linear
operators on a right quaternionic vector space.
\pn
\begin{definition} Let $V$ be a right vector space on $\mathbb{H}$.
A map $T: V\to V$ is said to be right linear if
$$
 T(u+v)=T(u)+T(v)
$$
$$
T(us)=T(u)s
$$
for all  $s\in\mathbb{H}$ and  all $u,v\in V$.
\end{definition}
\begin{remark} {\rm
Note that the set of right linear maps
is not a quaternionic left or right vector space.
Only if $V$ is both left and right vector space, then the set End$(V)$ of right linear maps on
$V$ is both a left and a right vector space on $\hh$, since in that case
we can define
$(aT)(v):=aT(v)$ and
$(Ta)(v):=T(av)$.
\noindent
The composition of operators can be defined in the usual way:
for any two operators $T,S\in {\rm End}(V)$ we have
$$
(TS)(u)=T(S(u)),\qquad \forall u\in V.
$$
In particular, we have the identity operator
$\mathcal{I}(u)=u$, for all $u\in V$ and
setting $T^0=\mathcal{I}$ we can define powers of a given operator $T\in {\rm End}(V)$:
$T^n=T(T^{n-1})$ for any $n\in \mathbb{N}$. An operator $T$ is said to be invertible
if there exists $S$ such that $TS=ST=\mathcal{I}$ and we will write $S=T^{-1}$.}
\end{remark}
\begin{remark}{\rm
From now on we will only consider bilateral vector spaces $V$.
The vector space End$(V)$ is not an $\hh$-algebra with respect to the composition of operators, in fact
the property $s(TS)=(sT)S=T(sS)$  is not fulfilled for any $T,S\in{\rm End}(V)$ and any
$s\in\hh$.
In this case we have $T(sS)(u)=T(sS(u))$ while $(sT)S(u))=sT(S(u))$.
Note on the other hand that End$(V)$ is trivially an algebra over $\rr$.
}
\end{remark}

\begin{definition} Let $V$ be a bilateral quaternionic Banach space.
We will denote by $\mathcal{B}(V)$ the bilateral vector space of all right linear
bounded operators on $V$.
\end{definition}
\noindent
It is easy to verify that
$\mathcal{B}(V)$ is a Banach space endowed with its natural norm.
\begin{definition}  An element $T\in \mathcal{B}(V))$ is said to be invertible if there exists
$T'\in \mathcal{B}(V)$ such that $TT'=T'T=\mathcal{I}$.
\end{definition}
\noindent
It is obvious that the set of all invertible elements of $\mathcal{B}(V)$
 is a group with respect to the composition of operators
defined in $\mathcal{B}(V)$.

\section{Algebraic properties of the Cauchy kernel}

This section is motivated by the fact that the Cauchy kernel can be defined
for non commuting variables $s$ and $q$ and
it is the key ingredient to define a functional calculus.
\begin{definition}\label{Cauchykernel}
Let $q$, $s\in \mathbb{H}$ such that $sq\not=qs$.
We will call non commutative Cauchy kernel series (shortly Cauchy kernel series)
the expansion
$$
S^{-1}(s,q):=\sum_{n\geq 0}q^ns^{-1-n},
$$
for $|q|<|s|$.
\end{definition}

\begin{theorem}\label{equazione per S}
Let $q$
and $s$ be two quaternions such that $qs\not=sq$ and consider
$$
S^{-1}(s,q):=\sum_{n\geq 0} q^n s^{-1-n}.
$$
Then the  inverse $S(s,q)$ of $S^{-1}(s,q)$ is solution to the equation
\begin{equation}\label{quadratica}
S^2+Sq-sS=0.
\end{equation}
\end{theorem}
\begin{proof}
Observe that
$$
S^{-1}(s,q)s=\sum_{n\geq 0} q^n s^{-1-n}s=\sum_{n\geq 0} q^n s^{-n}=1+qs^{-1}+q^2s^{-2}+\ldots
$$
and
$$
qS^{-1}(s,q)=q\sum_{n\geq 0} q^n s^{-1-n}=\sum_{n\geq 0} q^{1+n} s^{-1-n}=qs^{-1}+q^2s^{-2}+\ldots
$$
so that
$$
S^{-1}(s,q)s-qS^{-1}(s,q)=1
$$
keeping in mind that $S^{-1}S=SS^{-1}=1$ we get $$
S(S^{-1}s-qS^{-1})S=S^2
$$
from which we obtain the proof.
\end{proof}

\begin{lemma}
Let $R(s,q):=s-q$. Then $R(s,q)$ is a solution of equation (\ref{quadratica}) if and only if $sq=qs$.
\end{lemma}
\begin{proof}
This result follows immediately from the chain of equalities
$$
(s-q)^2+(s-q)q-s(s-q)=s^2-sq-qs+q^2-s^2+sq+sq-q^2=-qs+sq
$$
whose last term vanishes if and only if $sq=qs$.
\end{proof}
\begin{remark}
{\rm
If $s=s_0+s_1L$, $q=q_0+q_1L$ for some $L\in \mathbb{S}$, then
$sq=qs$.
}
\end{remark}

\begin{theorem} \label{thinvquat}
Let $q$, $s\in \mathbb{H}$ be such that
$qs\not=sq$. Then the non trivial solution of
\begin{equation}\label{Sequa}
S^2+Sq-sS=0
\end{equation}
 is given by
\begin{equation}\label{Ssolution}
S(s,q)=(s+q-2 \ Re [s])^{-1}(sq-|s|^2)-q=-(q-\overline s)^{-1}(q^2
-2q Re [s]+|s|^2).
\end{equation}
\end{theorem}

\begin{proof}
Let us begin by noticing that $S(s, q)$ is well defined: in fact, since $qs\ne sq$ by
hypothesis, we have $q\ne \overline s=-s+2\ Re [s]$ and hence
$(q-\overline s)=(s+q-2 \ Re [s])\ne 0$. Now, by Theorem \ref{equazione per S},
to prove the assertion it is enough to verify that $S(s,q)=(s+q-2 \ Re [s])^{-1}(sq-|s|^2)-q$
has $S^{-1}(s,q)$ as its inverse. Consider the equality
\begin{equation}\label{SS}
S(s,q)S^{-1}(s,q)=1
\end{equation}
which can be written as
$$
[(s+q-2 \ Re [s])^{-1}(sq-|s|^2)-q]\sum_{n\geq 0} q^n s^{-1-n}=1
$$
and also as
\begin{equation}\label{star}
(s+q-2 \ Re [s])^{-1}(sq-|s|^2)\sum_{n\geq 0} q^n s^{-1-n}-q\sum_{n\geq 0} q^n s^{-1-n}=1.
\end{equation}
If we multiply by $(s+q-2 \ Re [s])\ne 0$ both sides of (\ref{star}), we obtain the equivalent equality
$$
(sq-|s|^2)\sum_{n\geq 0} q^n s^{-1-n}-(s+q-2 \ Re [s])q\sum_{n\geq 0} q^n s^{-1-n}=s+q-2 \ Re [s]
$$
and
$$
sq\sum_{n\geq 0} q^n s^{-1-n}-|s|^2\sum_{n\geq 0} q^n s^{-1-n}
$$
$$
-sq\sum_{n\geq 0} q^n s^{-1-n}-q^2\sum_{n\geq 0} q^n s^{-1-n}+2 \ Re [s]q  \sum_{n\geq 0} q^n s^{-1-n}=s+q-2 \ Re [s].
$$
We therefore obtain
\begin{equation}\label{SS1}
(-|s|^2-q^2+2qRe [s])\sum_{n\geq 0} q^n s^{-1-n}=s+q-2 \ Re [s].
\end{equation}
Observing that $-|s|^2-q^2+2q Re [s]$ commutes with $q^n$ we can rewrite
this last equation as
\begin{equation}\label{fondamentale}
\sum_{n\geq 0} q^n(-|s|^2-q^2+2qRe [s]) s^{-1-n}=s+q-2 \ Re [s].
\end{equation}
Now the left hand side can be written as
$$
\sum_{n\geq 0} q^n(-|s|^2-q^2+2qRe [s]) s^{-1-n}=
$$
$$
(-|s|^2-q^2+2qRe [s]) s^{-1}
+q^1(-|s|^2-q^2+2qRe [s]) s^{-2}
+q^2(-|s|^2-q^2+2qRe [s]) s^{-3}
+...
$$
$$
=-\Big(|s|^2s^{-1}+q(-2s Re [s] +|s|^2)s^{-2}
+q^2(s^2-2s Re [s]+|s|^2)s^{-3}
+q^3(s^2-2s Re [s]+|s|^2)s^{-4}+...\Big).
$$
Using the identity
$$
s^2-2s Re [s]+|s|^2=0
$$
we get
$$
\sum_{n\geq 0} q^n(-|s|^2-q^2+2qRe [s]) s^{-1-n}=-|s|^2s^{-1}+qs^2s^{-2}=-|s|^2s^{-1}+q
$$
$$
=
-\overline{s}ss^{-1}+q=-\overline{s}+q=s-2 \ Re [s]+q
$$
which equals the right hand side of (\ref{SS1}), thus showing that (\ref{SS}) is an identity.

\noindent We can now pass to verify that
$$
S^{-1}(s,q)S(s,q)=1,
$$
which is
$$
\sum_{n\geq 0} q^n s^{-1-n}[(s+q-2 \ Re [s])^{-1}(sq-|s|^2)-q]=1,
$$
or equivalently
$$
\sum_{n\geq 0} q^n s^{-1-n} (s+q-2 \ Re [s])^{-1}(q^2 -2q Re [s]+|s|^2)=-1,
$$
$$
\sum_{n\geq 0} q^n s^{-1-n} (s+q-2 \ Re [s])^{-1}=-(q^2 -2q Re [s]+|s|^2)^{-1},
$$
$$
(q^2 -2q Re [s]+|s|^2)\sum_{n\geq 0} q^n s^{-1-n} =-(s+q-2 \ Re [s]),
$$
$$
\sum_{n\geq 0} q^n (q^2 -2q Re[s]+|s|^2)s^{-1-n} =-(s+q-2 \ Re [s]).
$$
The conclusion of the proof is reduced to verify (\ref{fondamentale}), that we have already done.
\end{proof}
\begin{remark}{\rm It is worth noticing that the proof of this last theorem does not rely on the
fact that the (real) components of $q$ commute. In fact, the theorem would hold even if
$q$ were to be in $\hh\otimes \mathbb{A}$ with $\mathbb{A}$ any non commutative algebra, for
example an algebra of
matrices on $\hh$. This fact will be exploited in Theorem \ref{thsmeno} where the variable $q$ is formally
replaced by an operator $T$ whose components do not necessarily commute. }
\end{remark}
\begin{remark}{\rm
Equation (\ref{Ssolution}) points out  that $S(s, q)=-(q-\overline s)^{-1}(q^2
-2q Re [s]+|s|^2)$ will have no inverse if  $(q^2 -2q Re [s]+|s|^2)= 0$. Set $s=s_0 +s_1I_s$, $(s_0, s_1 \in \rr)$. Since
\begin{equation}\label{luogo di zeri}
\{q\in \hh : (q^2 -2q Re [s]+|s|^2)=0\}=s_0+s_1\mathbb{S}
\end{equation}
then we conclude that $S(s, q)$ will have no inverse if $q\in s_0+s_1\mathbb{S}$.
A particular case appears when $q=\overline s$, and one could hope to be able
to split in factors the polynomial $(q^2 -2q Re [s]+|s|^2)$ and extend the inverse
of $S(s, q)$ up to the case $q=\overline s$.
However, as shown in the next results, this is not possible and the
function $S^{-1}(s,q)$ cannot be extended to a continuous function in $q=\overline{s}$.}
\end{remark}
\begin{lemma}\label{fattorizzazione} Let  $s\in \hh$ be a non real quaternion. Then there exists no degree-one quaternionic polynomial $Q(q)$ such that
\begin{equation}\label{radici}
q^2-2qRe[s]+|s|^2=(q-\overline s)Q(q).
\end{equation}
\end{lemma}
\begin{proof}
If such a polynomial $Q(q)$ exists, then the right-hand side of (\ref{radici})
would have a finite number of roots (namely $1$ or $2$), while the left-hand
side has infinitely many roots in view of (\ref{luogo di zeri}).
\end{proof}

\begin{theorem}
Let $q, s \in \mathbb{H}$. Then the function
$$
S^{-1}(s,q)=-(q^2-2qRe[s]+|s|^2)^{-1}(q-\overline{s})
$$
cannot be extended continuously to any point of the set
$$
\{ (s, q)\in \hh \times \hh \ :\ (q^2-2qRe[s]+|s|^2)=0\ \}.
$$
In particular, if  $qs\not= sq$ the limit
\begin{equation}\label{limite}
\lim_{q\to \overline{s}}S^{-1}(s,q)
\end{equation}
does not exist.
\end{theorem}
\begin{proof}
We prove that the limit (\ref{limite}) does not exist. Let $\varepsilon \in \mathbb{H}$
and consider
$$
S^{-1}(s,\overline{s}+\varepsilon)=
((\overline{s}+\varepsilon)^2-2(\overline{s}+\varepsilon)Re[s]+|s|^2)^{-1}\varepsilon
$$
$$
=((\overline{s}+\varepsilon)^2-2(\overline{s}+\varepsilon)Re[s]+|s|^2)^{-1}\varepsilon
=(\overline{s}\varepsilon+\varepsilon\overline{s} +\varepsilon^2
-2\varepsilon Re[s])^{-1}\varepsilon
$$
$$
=(\varepsilon^{-1}(\overline{s}\varepsilon+\varepsilon\overline{s} +\varepsilon^2
-2\varepsilon Re[s]))^{-1}
=(\varepsilon^{-1}\overline{s}\varepsilon+\overline{s} +\varepsilon
-2 Re[s]))^{-1}.
$$
If we now let $\varepsilon\to 0$, we obtain that the term
$\varepsilon^{-1}\overline{s}\varepsilon$ does not have a limit because
$$
\varepsilon^{-1}\overline{s}\varepsilon=
\frac{\overline\varepsilon}{|\varepsilon|^2}\overline{s}\varepsilon
$$
contains addends of type $\displaystyle\frac{\varepsilon_i\varepsilon_j s_\ell}{|\varepsilon|^2}$ with
$i,j,\ell\in \{0,1,2,3\}$ that do not have limit.
\end{proof}
We now define the non commutative Cauchy kernel and, with an abuse of notation, we will still denote it
by $S^{-1}(s,q)$.
\begin{definition}\label{Cauchykernel}
Let $q$, $s\in \mathbb{H}$ such that $sq\not=qs$.
We will call non commutative Cauchy kernel  (shortly Cauchy kernel)
the expression
$$
S^{-1}(s,q):=-(q^2-2qRe[s]+|s|^2)^{-1}(q-\overline{s}).
$$
\end{definition}

\section{The $S$-resolvent operator and the $S$-spectrum}
Let $T$ be a linear quaternionic operator. It is obvious that there are two natural eigenvalue
problems associated to $T$. The first, which one could call the left eigenvalue problem consists
in the solution of equation $T(v)=\lambda v$, and the second, which is called right eigenvalue problem,
and consists in the solution of the equation $T(v)=v\lambda$. We will discuss the associated
spectra later on, but the key observation is that none of them is
 useful to define a functional calculus. In this section we will identify the correct
 framework for the study of eigenvalues for quaternionic operators.
\begin{definition}(The $S$-resolvent operator series)
Let $T\in\mathcal{B}(V)$ and let $s\in \mathbb{H}$.
We define the $S$-resolvent operator series as
\begin{equation}\label{Sresolv}
S^{-1}(s,T):=\sum_{n\geq 0} T^n s^{-1-n}
\end{equation}
for $\|T\|< |s|$.
\end{definition}

\noindent If we denote by $\mathcal{I}$ the identity operator, we can state the following:

\begin{theorem}\label{thsmeno}
Let $T\in\mathcal{B}(V)$ and let $s \in \mathbb{H}$. Assume that
$\overline s$ is such that $T-\overline{s}\mathcal{I}$ is invertible. Then
\begin{equation}\label{STeq}
S(s,T)=(T-\overline s\mathcal{I})^{-1}(sT-|s|^2\mathcal{I})-T
\end{equation}
is the inverse of
$$
S^{-1}(s,T)=\sum_{n\geq 0} T^n s^{-1-n}.
$$
Moreover, we have
\begin{equation}\label{ciao}
\sum_{n\geq 0} T^n s^{-1-n}=-(T^2-2Re[s] T+|s|^2\mathcal{I})^{-1}(T-\overline{s}\mathcal{I}),
\end{equation}
for $\|T\|< |s|$.
\end{theorem}
\begin{proof}
The subsequent calculations have to be intended in the norm of the bounded
linear quaternionic operators,
the proof mimics the one of Theorem \ref{thinvquat}.
We verify that
\begin{equation}\label{esseessemenouno}
S(s,T)S^{-1}(s,T)=\mathcal{I}
\end{equation}
when $S(s,T)$ is given by (\ref{STeq}). Namely (\ref{esseessemenouno}) can be written as
$$
[(T+(s-2 \ Re [s])\mathcal{I})^{-1}(sT-|s|^2\mathcal{I})-T]\sum_{n\geq 0} T^n s^{-1-n}=\mathcal{I}
$$
and also as
$$
(T+(s-2 \ Re [s])\mathcal{I})^{-1}(sT-|s|^2\mathcal{I})\sum_{n\geq 0} T^n s^{-1-n}-
T\sum_{n\geq 0} T^n s^{-1-n}=\mathcal{I}.
$$
By applying  $T+(s-2 \ Re [s])\mathcal{I}$ to both hands sides, we obtain the equality
$$
(sT-|s|^2\mathcal{I})\sum_{n\geq 0} T^n s^{-1-n}-(T+(s-2 \ Re [s])\mathcal{I})
T\sum_{n\geq 0} T^n s^{-1-n}=T+(s-2 \ Re [s])\mathcal{I}
$$
which can be written as
$$
sT\sum_{n\geq 0} T^n s^{-1-n}-|s|^2\sum_{n\geq 0} T^n s^{-1-n}
$$
$$
-sT\sum_{n\geq 0} T^n s^{-1-n}-T^2\sum_{n\geq 0} T^n s^{-1-n}+2 \ Re [s]T
  \sum_{n\geq 0} T^n s^{-1-n}=T+(s-2 \ Re [s])\mathcal{I}
$$
and then we get
$$
(-|s|^2\mathcal{I}-T^2+2Re[s] T)\sum_{n\geq 0} T^n s^{-1-n}=T+(s-2 \ Re [s])\mathcal{I}.
$$
Observing that $-|s|^2\mathcal{I}-T^2+2Re[s] T$ commutes with $T^n$ we obtain that the above identity is equivalent to
$$
\sum_{n\geq 0} T^n(-|s|^2-T^2+2Re[s] T) s^{-1-n}=T+(s-2 \ Re [s])\mathcal{I}.
$$
Now we expand the series as
$$
\sum_{n\geq 0} T^n(-|s|^2\mathcal{I}-T^2+2Re[s] T) s^{-1-n}=
$$
$$
(-|s|^2\mathcal{I}-T^2+2 Re[s]T  ) s^{-1}
+T^1(-|s|^2\mathcal{I}-T^2+2 Re[s]T) s^{-2}
+T^2(-|s|^2\mathcal{I}-T^2+2 Re[s]T) s^{-3}
+...
$$
$$
=-\Big(|s|^2s^{-1}+T(-2s Re[s] +|s|^2)s^{-2}
+T^2(s^2-2s Re[s]+|s|^2)s^{-3}
+...\Big)
$$
and using the identity
$$
s^2-2s Re[s]+|s|^2=0
$$
we get
$$
\sum_{n\geq 0} T^n(-|s|^2-T^2+2Re[s]T) s^{-1-n}
=-|s|^2s^{-1}\mathcal{I}+Ts^2s^{-2}=-|s|^2s^{-1}\mathcal{I}+T
$$
$$
=
-\overline{s}ss^{-1}\mathcal{I}+T=-\overline{s}\mathcal{I}+T=(s-2 \ Re [s])\mathcal{I}+T.
$$
The equality (\ref{ciao}) follows directly by taking the inverse of (\ref{STeq}).
\end{proof}
\begin{theorem}\label{invalg}
Let $T\in\mathcal{B}(V)$ and let $s\in \mathbb{H}$.
Then the operator
$$
\sum_{n\geq 0}(s^{-1}T)^ns^{-1}\mathcal{I}
$$
is the right and left algebraic inverse of $s\mathcal{I}-T$. Moreover, the series converges
in the operator norm for $\|T\|< |s|$.
\end{theorem}
\begin{proof} Let us directly compute
$$
(s\mathcal{I}-T)\sum_{n\geq 0}(s^{-1}T)^ns^{-1}\mathcal{I}=
s\mathcal{I}\sum_{n\geq 0}(s^{-1}T)^ns^{-1}\mathcal{I}
  -T\sum_{n\geq 0}(s^{-1}T)^ns^{-1}\mathcal{I}
$$
$$
=
s\mathcal{I}s^{-1}\mathcal{I}  +Ts^{-1}\mathcal{I}
+T(s^{-1}T)s^{-1}\mathcal{I}  +T(s^{-1}T)^2s^{-1}\mathcal{I}+ \ldots
$$
$$
-Ts^{-1}\mathcal{I}-T(s^{-1}T)s^{-1}\mathcal{I}
-T(s^{-1}T)^2s^{-1}\mathcal{I}-T(s^{-1}T)^3s^{-1}\mathcal{I}+\ldots =\mathcal{I}.
$$
The same identity holds for
$$
\sum_{n\geq 0}(s^{-1}T)^ns^{-1}\mathcal{I}(s\mathcal{I}-T)=\mathcal{I}.
$$
Finally we consider
$$
\|\sum_{n\geq 0}(s^{-1}T)^ns^{-1}\mathcal{I}\|
\leq \sum_{n\geq 0}\|(s^{-1}T)^ns^{-1}\mathcal{I}\|\leq
\sum_{n\geq 0}\|(s^{-1}T)\|^n|s^{-1}|\leq \sum_{n\geq 0}\|T\|^n|s^{-1}|^{n+1}
$$
which converges for $\|T\|< |s|$.
\end{proof}

\begin{corollary}\label{44}
When $Ts\mathcal{I}=sT$, the operator $S^{-1}(s,T)$ equals $(s\mathcal{I}-T)^{-1}$
when the series (\ref{Sresolv}) converges.
\end{corollary}
\begin{proof} It follows immediately from Theorem \ref{invalg}.
\end{proof}
\begin{definition}(The $S$-resolvent operator)
Let $T\in\mathcal{B}(V)$ and let $s\in \mathbb{H}$.
We define the $S$-resolvent operator as
\begin{equation}\label{Sresolvoperator}
S^{-1}(s,T):=-(T^2-2Re[s] T+|s|^2\mathcal{I})^{-1}(T-\overline{s}\mathcal{I}).
\end{equation}
\end{definition}

\begin{definition}(The spectra of  quaternionic operators)
Let $T:V\to V$ be a linear quaternionic operator on the Banach space $V$.
We denote by $\sigma_L(T)$ the left spectrum of $T$ related to the resolvent
$(s\mathcal{I}-T)^{-1}$
that is
$$
\sigma_L(T)=\{ s\in \mathbb{H}\ \ :\ \ s\mathcal{I}-T\ \ \ {\it is\ not\  invertible} \}.
$$
We define the $S$-spectrum $\sigma_S(T)$ of $T$ related to the $S$-resolvent
 operator (\ref{Sresolvoperator}) as:
$$
\sigma_S(T)=\{ s\in \mathbb{H}\ \ :\ \ T^2-2 \ Re [s]T+|s|^2\mathcal{I}\ \ \
{\it is\ not\  invertible}\}.
$$
\end{definition}
\begin{remark}
{\rm
It is also possible to introduce a notion of right spectrum
$\sigma_R(T)$ of $T$ as
$$
\sigma_R(T)=\{ s\in \mathbb{H}\ \ :\ \ \mathcal{I}\cdot s-T\ \ \ {\rm is\ not\  invertible} \},
$$
where the notation $\mathcal{I}\cdot s$ means that the multiplication by $s$ is on the right, i.e.
$\mathcal{I}\cdot s(v)=\mathcal{I}(v) s$. However, the operator $\mathcal{I}\cdot
s-T$ is not linear, so we will never refer to this notion.
}
\end{remark}

\begin{theorem}
Let $T\in\mathcal{B}(V)$ and let $s \in \rho_S(T)$. Then the $S$-resolvent
operator defined in (\ref{Sresolvoperator}) satisfies the equation
$$
S^{-1}(s,T)s-TS^{-1}(s,T)=\mathcal{I}.
$$
\end{theorem}
\begin{proof}
It follows by direct computation. Indeed, replacing (\ref{Sresolvoperator})
in the above equation we have
\begin{equation}\label{serve21}
-(T^2-2Re[s] T+|s|^2\mathcal{I})^{-1}(T-\overline{s}\mathcal{I})s
+T
(T^2-2Re[s] T+|s|^2\mathcal{I})^{-1}(T-\overline{s}\mathcal{I})=\mathcal{I}
\end{equation}
and applying $(T^2-2Re[s] T+|s|^2\mathcal{I})$ to both hands sides of (\ref{serve21}),
we get
$$
-(T-\overline{s}\mathcal{I})s+(T^2-2Re[s] T+|s|^2\mathcal{I})T
(T^2-2Re[s] T+|s|^2\mathcal{I})^{-1}(T-\overline{s}\mathcal{I})=T^2-2Re[s] T+|s|^2\mathcal{I}.
$$
Since $T$ and $T^2-2Re[s] T+|s|^2\mathcal{I}$ commute, we obtain the identity
$$
-(T-\overline{s}\mathcal{I})s+T
(T-\overline{s}\mathcal{I})=T^2-2Re[s] T+|s|^2\mathcal{I}
$$
which proves the statement.
\end{proof}
\begin{definition} The equation
$$
S^{-1}(s,T)s-TS^{-1}(s,T)=\mathcal{I}
$$
will be called the $S$-resolvent equation.
\end{definition}
\begin{theorem}
Let $T\in\mathcal{B}(V)$.
Then $\sigma_S(T)$ and $\sigma_L(T)$ are contained in the set
$
\{  s\in \mathbb{H}  : |s|\leq \|T\| \}.
$
\end{theorem}
\begin{proof}
Since both the series
$$
\sum_{n\geq 0}(s^{-1}T)^ns^{-1}\mathcal{I},\ \ \ \sum_{n\geq 0} T^n s^{-1-n}
$$
converge if and only if $|s^{-1}| \|T\|<1$ we get the  statement.
\end{proof}
\begin{theorem} Let $T\in\mathcal{B}(V)$ and $s\not\in \sigma_S(T)$ such that
\begin{equation}\label{stiaa}
\|T-Re [s] \mathcal{I} \|<|s-Re[s]|.
\end{equation}
Then the $S$-resolvent operator admits the series expansion:
\begin{equation}\label{stiaaaa}
S^{-1}(s,T)=\sum_{n\geq 0} (T-Re [s] \mathcal{I})^{n}
(s-Re[s])^{-n-1}.
\end{equation}
\end{theorem}
\begin{proof}
Observe that
$$
(T^2-2TRe[s]+|s|^2\mathcal{I})^{-1}=[(T-Re[s]\mathcal{I})^2+|s|^2\mathcal{I}-(Re[s])^2\mathcal{I}]^{-1}
$$
\begin{equation}\label{sti1}
=(|s|^2-(Re[s])^2)^{-1}\Big[\mathcal{I}+\frac{(T-Re[s])^2}{(|s|^2-(Re[s])^2)}\Big]^{-1}
=\sum_{n\geq
0}(-1)^n\frac{(T-Re[s])^{2n}}{(|s|^2-(Re[s])^2)^{n+1}}.
\end{equation}
Since $|Im [s]|^{2}=-(Im [s])^{2}$, then by replacing
(\ref{sti1}) in (\ref{Sresolvoperator}) we get
$$
S^{-1}(s,T)=\sum_{n\geq 0} (T-Re[s]\mathcal{I})^{2n+1}
(s-Re[s])^{-2n-2}+ \sum_{n\geq 0}(T-Re[s]\mathcal{I})^{2n}
(s-Re[s])^{-2n-1}
$$
adding the two terms we get (\ref{stiaaaa}) which converges when
(\ref{stiaa}) holds.
\end{proof}

We now give a simple relation between the $L$-spectrum and the $S$- spectrum.
\begin{prop}
Let $T\in\mathcal{B}(V)$ and $s\in \sigma_L(T)$ and let $v$ be the corresponding  $L$-eigenvector.
Then $s\in \sigma_S(T)$ and $v$ is the corresponding $S$-eigenvector if and only if
$$
(T- s\mathcal{I}) (s v)=0.
$$
\end{prop}
\begin{proof} It follows from the relations:
$$
T^2v-2 Re[s] Tv+\overline{s} s\mathcal{I} v=
T(sv)-2 Re[s] (sv)+\overline{s} (s v)=
(T- s\mathcal{I}) (s v)=0.
$$
\end{proof}
\noindent
We now announce a key result, whose proof we must postpone to Section 5.
\begin{theorem}\label{compattezaS}(Compactness of $S$-spectrum)
 Let $T\in\mathcal{B}(V)$. Then
the $S$-spectrum $\sigma_S (T)$  is a compact nonempty set contained in
$\{\ s\in \mathbb{H}\ :\  |s|\leq \|T\| \ \}$.
\end{theorem}
\begin{theorem}\label{strutturaS} (Structure of the $S$-spectrum)
 Let $T\in\mathcal{B}(V)$ and let $p = p_0 +p_1I\in   p_0 +p_1\mathbb{S}\subset \hh\setminus\rr$
 be an $S$-eigenvalue of $T$.
  Then all the elements of the sphere $p_0 +p_1\mathbb{S}$
 are $S$-eigenvalues of $T$.
\end{theorem}
\begin{proof}
In the $S$-eigenvalue equation the coefficients depend only on the
real numbers $p_0, p_1$ and not on $I \in \mathbb{S}$. Therefore all $s=p_0+p_1J$ such that $J \in \mathbb{S}$
 are in the $S$-spectrum.
\end{proof}
\noindent
We conclude this section with a few examples which highlight the differences between
the $S$-spectrum and the $L$-spectrum.
\begin{example}{\rm
1. Consider the matrix
$$
T_1=\left[%
\begin{array}{cc}
  1 &0  \\
0 & j%
\end{array}
\right].
$$
Direct computations show that $\sigma_S(T_1)=\{1\}\cup \mathbb{S}$,
while  $\sigma_L(T_1)=\{1, j\}$.
\pn
2.
We now consider
$$
T_2=\left[%
\begin{array}{cc}
  i &0  \\
0 & j%
\end{array}
\right].
$$
We have (as a point-set) $\sigma_S(T_2)=\mathbb{S}$
and, obviously, the $L$-spectrum is $\sigma_L(T_2)=\{ \ i,\ \ j \ \}.$\pn
3. Finally, let
$$
T_3=\left[%
\begin{array}{cc}
  0 &i  \\
-i & 0
\end{array}
\right].
$$
We obtain
$
\sigma_S(T_3)=\{ \pm 1   \} $ and $\sigma_L(T_3)=\{ s_0+s_1i+s_2j+s_3k :s_1=0,\ \ \ s_0^2+s_2^2+s_3^2=1\} = \{s\in \mathbb{H} : si\in \mathbb{S}\}.
$}
\end{example}

\section{Functional calculus}
\begin{definition} A function  $f:\ \hh \to \hh$
is said to be locally
regular on the spectral set $\sigma_S(T)$ of an operator $T\in \mathcal{B}(V)$
if there exists a ball $ B(0,R)$ containing $\sigma_S(T)$
on which $f$ is regular.
We will denote by ${\cal R}_{\sigma_S(T)}$
the set of locally regular functions on $\sigma_S (T)$.
\end{definition}
\begin{theorem} Let $T\in\mathcal{B}(V)$ and  $f\in {\cal R}_{\sigma_S(T)}$. Choose $I\in \mathbb{S}$ and
let $L_I$ be the plane that contains the real line and the
imaginary unit $I$. Let $U$ be an open bounded set in $L_I$ that
contains $L_{I}\cap \sigma_S(T)$. Then the integral
$$
{{1}\over{2\pi }} \int_{\partial U } S^{-1} (s,T)\  ds_I \ f(s)
$$
does not depend on the choice of the imaginary unit $I$ and on the open set $U$.
\end{theorem}
\begin{proof}
By Theorem \ref{strutturaS}
the $S$-spectrum contains either real points or (entire) spheres of  type
$s_0+r\mathbb{S}$, $(s_0, r \in \rr)$.
Every plane $L_J=\mathbb{R}+J\mathbb{R}$, $(J\in \mathbb{S})$, contains all the real points of the $S$-spectrum. Moreover, $L_J$ intersects each sphere $s_0+r\mathbb{S}$ in the two (conjugate) points $s_0+rJ$ and $s_0-rJ$.
Now we show that the integral in the statement does not
depend on the plane $L_I$ and on $U$.
Now,
for any two imaginary units $I, I'\in\mathbb{S}$,  $I'\neq I$, and any $U'$ containing
$L_{I'}\cap \sigma_S(T)$ we obtain:
$$
{{1}\over{2\pi }}\int_{\partial U} \sum_{n\geq 0} T^n s^{-1-n}\  ds_I f(s)=
{{1}\over{2\pi }}\int_{\partial U'} \sum_{n\geq 0} T^n s^{-1-n}\  ds_{I'} f(s)
$$
thanks the the Cauchy theorem and thanks to the fact that the points of the spectrum
have that same coordinates on each ``complex" plane $L_J$.
\end{proof}
We give a preliminary result that motivates the functional calculus.
\begin{theorem}\label{poly}
Let $ q^m a$ be a monomial, $q, a\in \mathbb{H}$, $m\in \mathbb{N}\cup \{0\}$.
Let $T\in\mathcal{B}(V)$ and let $U$ be an open bounded set
  in $L_I$ that contains  $L_{I}\cap \sigma_S(T)$.
 Then
\begin{equation}\label{TA}
T^m a= {{1}\over{2\pi }} \int_{\partial U } S^{-1} (s,T)\  ds_I \  s^m \ a.
\end{equation}
\end{theorem}
\begin{proof}
Let us consider the power series expansion for the operator $S^{-1} (s,T)$ and a circle $C_r$
in $L_I$
centered at 0 and with radius $r>\|T\|$. We have
\begin{equation}\label{25}
{{1}\over{2\pi }}\int_{C_r} S^{-1} (s,T)\  ds_I \  s^m \ a
={{1}\over{2\pi }}\sum_{n\geq 0} T^n\int_{C_r} s^{-1-n+m}\  ds_I a
 = T^m \ a,
\end{equation}
since the following equalities hold:
\begin{equation}
\int_{C_r} ds_I  s^{-n-1+m}=0\ \ if \ \ n\not=m, \ \ \ \ \ \
\int_{C_r} ds_I  s^{-n-1+m}=2\pi \ \ if \ \ n=m.
\end{equation}
The Cauchy theorem shows that the integral (\ref{25}) is not affected if  $C_r$ is replaced by $\pp U$.
\end{proof}
\noindent
The preceding discussion allows to give the following definition.
\begin{definition} \label{fdiT}
Let $f$ be a regular function and $T$ a linear bounded operator.
We define
$$
f(T)= {{1}\over{2\pi }} \int_{\partial U } S^{-1} (s,T)\  ds_I \ f(s)
$$
where $U$ is an open bounded set that contains $L_{I}\cap \sigma_S(T)$.
\end{definition}
\noindent
We are now ready to prove the compactness of the spectrum
which we announced in the previous section.\par\noindent
{\bf Theorem \ref{compattezaS}} {\it (Compactness of $S$-spectrum)
 Let $T\in\mathcal{B}(V)$. Then
the $S$-spectrum $\sigma_S (T)$  is a compact nonempty set
contained in $\{\ s\in \mathbb{H}\ :\  |s|\leq \|T\| \ \}$.
}
\begin{proof}
Observe that
$$
\frac{1}{2\pi }\int_{\partial U} \sum_{n\geq 0} T^n s^{-1-n}\
ds_I \  s^m\, =T^m
$$
where $U$ contains the part of the $S$-spectrum in the plane $L_I$.
So for $m=0$
$$
\frac{1}{2\pi }\int_{\partial U } \sum_{n\geq 0} T^n s^{-1-n}\
ds_I \  =\mathcal{I},
$$
where $\mathcal{I}$ denotes the identity operator. This means that
$\sigma_S (T) $ is a non empty set. We show that it is bounded.
The series $ \sum_{n\geq 0} T^n s^{-1-n} $ converges if and only
if $\|T\|< |s|$ so the $S$-spectrum is contained in the set $\{s
\in \mathbb{H}\ :\ |s| \leq \|T\| \}$, which is bounded and closed
because
 the complement set of $\sigma_S (T) $
is open. Indeed, the function $g: s\mapsto T^2-2 Re[s] T+|s|^2\mathcal{I}$
is trivially continuous and, by Theorem 10.12 in \cite{rudin}, the set  $\mathcal{U}(V)$
of all invertible elements of $\mathcal{B}(V)$  is an open set in $\mathcal{B}(V)$.
Therefore $g^{-1}(\mathcal{U}(V))=\rho_S(T)$ is an open set in $\mathbb{H}$.

\end{proof}
The next two theorems show that Definition \ref{fdiT} leads to a functional calculus with good
properties.
\begin{theorem}
Let $T\in\mathcal{B}(V)$ and let $f$
and $g\in {\cal R}_{\sigma_S(T)}$. Then
$$
(f+g)(T)=f(T)+g(T),\ \ \ \ \ (f\lambda )(T)=f(T) \lambda  \ \ \ \ \
{\it for\ all}\ \ \lambda  \in \mathbb{H}.
$$
Moreover, if $\phi(s)=\sum_{n\geq 0} s^n a_n$ and  $\psi(s)=\sum_{n\geq 0} s^n b_n$,
 are  in $ {\cal R}_{\sigma_S(T)}$ with $a_n$
 and $b_n\in \mathbb{R}$.
Then
$$
(\phi\psi)(T)=\phi(T)\psi(T).
$$
\end{theorem}
\begin{proof}
The first part of the theorem is a direct consequence of Definition \ref{fdiT}.
Observe that when $a_n$
 and $b_n$ are real the product $\phi\psi$ is regular.
We have to prove that
$$
\phi(T) \psi(T)=\frac{1}{2\pi }\int_{\partial U}S^{-1}(s,T) \, ds_I\,\phi(s)\psi(s),
$$
where $U$ contains $L_{I}\cap \sigma_S(T)$. By
recurrence for $m\geq 2$ we get
$$
S^{-1}(s,T)s^m- T^mS^{-1}(s,T)=
\mathcal{I}s^{m-1}+Ts^{m-2}+T^2s^{m-3}+...+T^{m-1},
$$
set
$$
Q_m(s,T):=\mathcal{I}s^{[m-1]_+}+Ts^{[m-2]_+}+T^2s^{[m-3]_+}+...+T^{[m-1]_+}, \ \ \ m\geq 2
$$
where $s^{[n]_+}=s^n$ if $n> 0$ , $s^{[n]_+}=0$ otherwise,
and
$$
Q_m(s,T):=\mathcal{I}\ \  \ \ \ {\rm for}\ \ \ m=1\ \ \ \ {\rm and}\ \ \ \ Q_m(s,T):=0\ \  \ \ \ for\ \ \ m=0,
$$
which is a regular function in $s$, with values in the space of linear bounded
quaternionic operators, so  we get
$$
T^mS^{-1}(s,T)=S^{-1}(s,T)s^m -Q_m(s,T),\ \ \ \ m=0,1,2,....
$$
We have to observe that for $m=1$ and $m=0$ we have
$ TS^{-1}(s,T)=S^{-1}(s,T)s-\mathcal{I}$ and $\mathcal{I}S^{-1}(s,T)=S^{-1}(s,T)1$, respectively.
Observe that
$$
\psi(T)=\frac{1}{2\pi }\int_{\partial U}S^{-1}(s,T) \, ds_I\, \psi(s)
$$
and
$$
T^m\psi(T)=T^m \frac{1}{2\pi }\int_{\partial U}S^{-1}(s,T) \, ds_I\, \psi(s)
=\frac{1}{2\pi }\int_{\partial U}T^mS^{-1}(s,T) \, ds_I\, \psi(s)
$$
$$
=\frac{1}{2\pi }\int_{\partial U}[S^{-1}(s,T)s^m -Q_m(s,T)] \, ds_I\, \psi(s)
=\frac{1}{2\pi }\int_{\partial U}S^{-1}(s,T)s^m  \, ds_I\, \psi(s)
$$
Now for $a_n\in \mathbb{R}$ we can write
$$
T^m a_m\psi(T)=\frac{1}{2\pi }
\int_{\partial U}S^{-1}(s,T)s^ma_m  \, ds_I\, \psi(s)
$$
summing up with respect to $m$
$$
\sum_{m=0}^MT^m a_m\psi(T)=
\frac{1}{2\pi }\int_{\partial U}S^{-1}(s,T)\sum_{m=0}^Ms^ma_m  \, ds_I\, \psi(s)
$$
now by definition $|\sum_{m=0}^Ms^ma_m|$ is bounded by $|\phi|$ which
is continuous on the bounded set $\partial U$,
and $\|S^{-1}(s,T)\|\leq (|s|-\|T\|)^{-1}$ so we can  pass
 to the limit for $M\to\infty$.

\end{proof}

\begin{theorem}(Spectral decomposition of a quaternionic  operator)
Let $T\in\mathcal{B}(V)$.
 Let $L_I\cap\sigma_S(T)= \sigma_{1S}(T)\cup \sigma_{2S}(T)$,
with $dist( \sigma_{1S}(T),\sigma_{2S}(T))>0$. Let $U_1$ and $U_2$
be two open sets such that  $\sigma_{1S}(T)
\subset U_1$ and $ \sigma_{2S}(T)\subset U_2$ , on $L_I$, with
$\overline{U}_1 \cap\overline{U}_2=\emptyset$. Set
$$
P_j:=\frac{1}{2\pi }\int_{\partial U_j}S^{-1}(s,T) \, ds_I\, \ \ \ \ \ \
T^m_j:=\frac{1}{2\pi }\int_{\partial U_j}S^{-1}(s,T) \, ds_I\, s^m\,
,\ \ \ m\in\mathbb{N},\ \ j=1,2.
$$
Then $P_j$ are projectors and
\begin{itemize}
\item[(I)] \ \ \  $P_1+P_2=\mathcal{I}$,

\item[(II)] \ \ \ $TP_j=T_j$,

\item[(III)] \ \ \ $T=T_1+T_2$,

\item[(IV)] \ \ \ $T^m=T^m_1+T^m_2$, \ $m \geq 2$.

\end{itemize}
\end{theorem}
\begin{proof}
Observe that $P_j=T_j^0$ and  note that the resolvent equation for $m=0$ is trivially
$T^0_jS^{-1}(s,T)=S^{-1}(s,T)s^0=S^{-1}(s,T)$. We have
$$
P_j^2=P_j\frac{1}{2\pi }\int_{\partial U_j}S^{-1}(s,T) \, ds_I\,
=
\frac{1}{2\pi }\int_{\partial U_j}P_jS^{-1}(s,T) \, ds_I\,
=\frac{1}{2\pi }\int_{\partial U_j}S^{-1}(s,T) \, ds_I\, =P_j.
$$
To prove (I) we use the Cauchy integral theorem, if $ \overline{U}_1 \cup\overline{U}_2 \subset U$, then
$$
\frac{1}{2\pi }\int_{\partial U_1}S^{-1}(s,T) \, ds_I\,
+\frac{1}{2\pi }\int_{\partial U_2}S^{-1}(s,T) \, ds_I \, =
\frac{1}{2\pi }\int_{\partial U}S^{-1}(s,T) \, ds_I .
$$
This gives $P_1+P_2=\mathcal{I}$ since
$\frac{1}{2\pi }\int_{\partial U}S^{-1}(s,T) \, ds_I =\mathcal{I}$.
\par\noindent
To prove (II) we  recall the resolvent relation
$$TS^{-1}(s,T)=S^{-1}(s,T)s-\mathcal{I}$$
so
$$
TP_j=\frac{1}{2\pi }\int_{\partial U_j} TS^{-1}(s,T) \, ds_I =
\frac{1}{2\pi }\int_{\partial U_j} [S^{-1}(s,T)s-\mathcal{I}] \, ds =
\frac{1}{2\pi }\int_{\partial U_j} S^{-1}(s,T) \, ds_I\, s\, =T_j.
$$
Now adding the relations $T_j=TP_j$ we get
$$
T_1+  T_2=TP_1+TP_2=T(P_1+P_2)=T,
$$
where we have used (I).
\par
\noindent
Now we observe that
$$
S^{-1}(s,T)s^2- T^2S^{-1}(s,T)=\mathcal{I}s+T
$$
and
$$
S^{-1}(s,T)s^3- T^3S^{-1}(s,T)=\mathcal{I}s^2+Ts+T^2
$$
and by recurrence for $m\geq 2$ we get
$$
S^{-1}(s,T)s^m- T^mS^{-1}(s,T)=
\mathcal{I}s^{m-1}+Ts^{m-2}+T^2s^{m-3}+...+T^{m-1}.
$$
Now, for $m\geq 2$, consider
$$
T^mP_j=\frac{1}{2\pi }\int_{\partial U_j} T^mS^{-1}(s,T) \, ds_I
$$
$$=
\frac{1}{2\pi }\int_{\partial U_j} [S^{-1}(s,T)s^m-(\mathcal{I}s^{m-1}
+Ts^{m-2}+T^2s^{m-3}+...+T^{m-1})] \, ds_I
$$
$$
=
\frac{1}{2\pi }\int_{\partial U_j} S^{-1}(s,T) \, ds_I\, s^m=T^m_j.
$$
So adding $T^mP_1=T^m_1$ and $T^mP_2=T^m_2$ and recalling (I) we get (IV).
\end{proof}


\begin{thebibliography}{99}

\bibitem{adler}
 S. Adler, {\it Quaternionic Quantum Field Theory}, Oxford University Press
(1995).


\bibitem{nostro libro} F. Colombo, I. Sabadini, F. Sommen, D.C.
Struppa, {\em Analysis of Dirac Systems and Computational Algebra},
Progress in Mathematical Physics, Vol. 39, {Birkh\"auser}, Boston,
2004.

\bibitem{ds}
N. Dunford, J. Schwartz, {\it Linear operators, part I: general
theory}, J. Wiley and Sons (1988).

\bibitem{gs} G. Gentili, D.C. Struppa, {\em A new approach to Cullen-regular functions
of a quaternionic variable},
 C.R. Acad. Sci. Paris, {\bf 342} (2006), 741--744.

\bibitem{advances} G. Gentili, D.C. Struppa, {\em A new theory of regular functions
of a quaternionic variable}, Advances in Mathematics, 2007, to
appear.

\bibitem{jefferies} B. Jefferies, {\em Spectral properties of noncommuting operators},
Lecture Notes in Mathematics, 1843, Springer-Verlag, Berlin, 2004.

\bibitem{rudin} W. Rudin, {\em Functional Analysis},
Functional analysis. McGraw-Hill Series in Higher Mathematics.
McGraw-Hill Book Co., New York-D\"usseldorf-Johannesburg, 1973.




\bibitem{Taylor} J.L. Taylor, {\em The analytic-functional calculus for several
 commuting operators},
Acta Math., {\bf 125} (1970), 1--38.
\bibitem{Taylor3} J.L. Taylor,
{\em Functions of several noncommuting variables}, Bull. Amer.
Math. Soc., {\bf 79} (1973), 1--34.


\end{thebibliography}
\end{document}